\documentclass[a4paper,twoside,10pt,reqno]{article}
\usepackage[english]{babel}
\usepackage[dvips]{graphicx}
\usepackage[T1]{fontenc}
\usepackage[latin1]{inputenc}
\usepackage{amssymb,amsthm,dsfont,mathrsfs}
\usepackage{amsfonts,amsmath,euscript}
\usepackage{color,cite}

\newtheorem{defi}{Definition}[section]
\newtheorem{teo}{Theorem}[section]
\newtheorem{pro}[teo]{Proposition}

\newtheorem{ex}{Example}[section]

\newcounter{example}[section]

           \addto{\captionsenglish}{}

\newcommand{\hs}{\hspace{3pt}}

\newcommand{\dem}{{\bf Proof. }}
\newcommand{\fdem}{$\square$}

\newcommand{\titulo}[1]{\mbox{} \\ \noindent \textit{\textbf{\Large #1}}\\}

\renewcommand{\abstract}[1]{{\small \noindent \textbf{Abstract:} #1\\}}

\newcommand{\pchave}[1]{{\small \noindent \textbf{Keywords:} #1\\}}

\begin{document}

\begin{center}
\titulo{Max-min dependence coefficients\\
for Multivariate Extreme Value Distributions}
\end{center}

\vspace{0.5cm}

\textbf{Helena Ferreira} Department of Mathematics, University of Beira Interior, Covilhã, Portugal\\

\abstract{We measure the dependence among sub vectors of a random vector with Multivariate Extreme Value distribution by using the expected value of a range and relate this coefficient of dependence with the multivariate tail dependence and extremal coefficients. The introduced coefficient extends the concept of madogram for several locations and several regions. The results are illustrated with some usual distributions and applied to financial data.}

\pchave{multivariate extreme value theory, dependence coefficients, range}

\section{Introduction}

The dependence structure of a Multivariate Extreme Value (MEV) distribution is completely characterised by its dependence function (Resnick (1987), Beirlant et al. (2004)). Since this function cannot be easily inferred from data the dependence coefficients are useful, despite the fact that one  coefficient cannot preserve all the information about this function.

\noindent The most popular of the dependence coefficients are those based on the tail dependence (Sybuya (1960), Li (2009)). They resume the probability of occurrence of extreme values for one or more random variable given that another(s) assumes  extreme values too. For the MEV distributions the extremal coefficient (Tiago de Oliveira (1962-63), Smith (1990)) is certainly a crucial and perhaps insurmountable tool when we have to summarize the dependence. For a $d-$dimensional random vector we have $2^d-d$ extremal coefficients which consistency pro\-per\-ties are discussed in Schlather and Tawn (2002). For an overview of other dependence measures for dependence see, for instance, Joe (1997).

\noindent  To the best of our knowledge, there is no extremal dependence coefficients for $p\geq 3$ subvectors  ${\bf X}_{1},...,{\bf X}_{p}$ of a random vector ${\bf X}=(X_1,...,X_d)$ with MEV distribution. 
 
\noindent  The need of to evaluate the strength of dependence among sub vectors arises for instance in the setting of max-stable random fields. Let $\{X_{{\bf{i}}}\}_{{\bf{i}}\in \mathbb{R}^2}$ be a max-stable random field and   $I_1,...I_p$ sets of locations in $\mathbb{R}^2$. The joint distribution of $X_{{\bf{i}}}, {\bf{i}}\in \displaystyle\bigcup_{j=1}^pI_j$, is a MEV distribution and we want to resume the dependence among the grouped values $\{X_{{\bf{i}}}, {\bf{i}}\in I_j\}$ at different regions $I_j$, $j=1,...,p$. This problem is treated by several authors for two variables $X_{{\bf{i}}}$ and  $X_{{\bf{j}}}$ corresponding to two locations ${\bf{i}}$  and ${\bf{j}}$ (see Naveau et al. (2009) and references therein) and the obtained results are extended for two regions $I_1$ and $I_2$ in Fonseca et al. (2012).
 
\noindent  In finance, we are frequently interested in assessing the dependence among several big world markets, considering each one as a random subvector. For an application with grouped financial stock markets see for instance  Ferreira and Ferreira (2012a).
 
\noindent  We propose to evaluate the degree of dependence among sub-vectors ${\bf X}_{1},...,{\bf X}_{p}$ of ${\bf X}$ with MEV distribution by using an expected range, which will be referred as a "max-min coefficient". It is a summary measure that takes into account the whole group of the extremal coefficients $\epsilon_{{\bf X}_{j}}$ of $X_{{\bf{j}}}$, $j=1,...,p$. Our approach is an extension of the modeling for pairwise dependence throughout the madogram (Poncet et al. (2006), Naveau et al. (2009)), an extreme-value analogue of the variogram (Cressie (1993)), since it enables to resume the spatial dependence structure for several locations or regions of locations. 

\noindent The proposed moment-based dependence tool takes into account the spread and dependence among the subvectors and can be easily estimated.

\noindent The paper is organized as follows. We introduce in section 2 the dependence coefficient whose is well defined for any random vector with MEV distribution, is a function of its copula and is invariant with respect permutations of the variables. Its relations with the multivariate tail dependence and the extremal coefficients are presented. Based on the expected range coefficient considered we compare a MEV distribution with others more concordant distributions and state some bounds.

\noindent In section 3, we  compute the max-min coefficients for the marginal distribution of the Multivariate Maxima of Moving Maxima process and the Symmetric Logistic distribution. We refer briefly an estimator for the max-min coefficients and apply it to grouped financial stock markets.

%postpone a comparative study with applications to a future work.
%__________________________________________________________
\section{Max-min dependence coefficients}\label{sed}

Let ${\bf X}=(X_1,...,X_d)$ be a  vector of unit Fr\'{e}chet random variables, that is, with marginal distribution function $F(x)=\exp (-x^{-1})$, $x>0$,
and $G$ denote the Multivariate Extreme Value distribution of ${\bf X}$.

\noindent The tail dependence function (Huang (1992), Schmidt and Stadtmuller (2006)) of $G$ is defined by $l(x_1,...,x_d)=-\log G(x_1^{-1},...,x_d^{-1})$, $(x_1,...,x_d) \in [0,\infty)^d$. It is a convex function, homogeneous of order one and satisfies
\begin{eqnarray*}\label{01}
\displaystyle\bigvee _{j=1}^d x_j\leq l(x_1,...,x_d)\leq \sum_{j=1}^dx_j,
\end{eqnarray*}

\noindent with the lower bound corresponding to ${\bf X}$ with totally dependent margins and the upper bound to ${\bf X}$ with independent margins.

\noindent The tail dependence function $l(x_1,...,x_d)$ gives us information about the probability of occurrence of extreme events for the maximum $\displaystyle\bigvee _{j=1}^dF(X_j)$ given that one fixed margin assumes extreme values too. In fact, we have (Schmidt and Stadmuller (2006), Ferreira and Ferreira(2011))

$$l(x_1,...,x_d)=\displaystyle\lim_{t\rightarrow\infty}\left(-t \log P\left(F(X_1)\leq 1-\frac{x_1}{t},...,F(X_d)\leq 1-\frac{x_d}{t}\right)\right)=$$
$$\displaystyle\lim_{t\rightarrow\infty} tP\left(F(X_1)>1- \frac{x_1}{t}\vee...\vee F(X_d)>1- \frac{x_d}{t}\right)$$

\noindent and, in particular,

$$l(x,...,x)=xl(1,...,1)=x\displaystyle\lim_{t\rightarrow\infty} P\left(\displaystyle\bigvee _{j=1}^dF(X_j)>1-\frac{1}{t}\mid F(X_i)>1-\frac{1}{t}\right).$$

\noindent For $\delta_i(S)=1$ if $i\in S$ and $\delta_i(S)=0$ if $i\notin S$, it holds that 
\begin{eqnarray}\label{2a}
l(\delta_1(S),...,\delta_d(S))=\epsilon_{{\bf X}_S}, 
\end{eqnarray}
\noindent where $\epsilon_{{\bf X}_S}$ is the extremal coefficient of the subvector ${\bf X}_S$ of ${\bf X}$ with indices in $S$ (Tiago de Oliveira (1962-63), Smith (1990)). It takes values in $[1,|S|]$, with $\epsilon_{{\bf X}_S}=1$ when ${\bf X}_S$ has the minimum copula $C_{{\bf X}_S}(u_1,...,u_d)_S=\displaystyle\bigwedge_{j\in S}u_j$ and $\epsilon_{{\bf X}_S}=|S|$ when ${\bf X}_S$ has the product copula $C_{{\bf X}_S}(u_1,...,u_d)_S=\displaystyle\prod_{j\in S}u_j$.

\noindent  Let $\mathcal{I}=\{I_1,...,I_p\}$ be a partition of $D=\{1,...,d\}$, $M(I_j)=\displaystyle\bigvee_{i\in I_j}X_i$ and  ${\bf X}_{I_j}$ the subvector of ${\bf X}$ with indices in $I_j$.
We resume the extremal dependence among ${\bf X}_{I_j}$, $j=1,...,p$, by the coefficient $R({\bf X},{\bf \lambda},\mathcal{I})$ defined as follows.\\
%-

 \begin{defi}Let ${\bf X}$ be a vector of unit Fr\'{e}chet random variables and Multivariate Extrem Value distribution. For each ${\bf \lambda}=(\lambda_1,...,\lambda_p)\in (0,\infty)^p$ and partition $\mathcal{I}$ of $D$, we define

$$R({\bf X},{\bf \lambda},\mathcal{I})=E\left(\displaystyle\bigvee _{j=1}^pF^{\lambda_j}(M(I_j))-\displaystyle\bigwedge _{j=1}^pF^{\lambda_j}(M(I_j))\right).$$
\end{defi}
\vspace{0.5cm}
\noindent  We remark the relations
\begin{eqnarray}\label{8a}
R({\bf X},{\bf \lambda},\mathcal{I})=E\left(\displaystyle\bigvee _{\{i,j\}\subset \{1,...,p\}}\left |F^{\lambda_j}(M(I_j))-F^{\lambda_i}(M(I_i))\right| \right)
\end{eqnarray}
\noindent  and 

\begin{eqnarray}\label{8b}
R({\bf X},{\bf \lambda},\mathcal{I})=E\left(\displaystyle\bigvee _{j=1}^p\displaystyle\bigvee _{i\in I_j}F^{\lambda_j}(X_i)-\displaystyle\bigwedge _{j=1}^p\displaystyle\bigvee _{i\in I_j}F^{\lambda_j}(X_i)\right).
\end{eqnarray}

\noindent By taking $p=2$ in (\ref{8a}), we find in $\frac{1}{2}R({\bf X},{\bf \lambda},\mathcal{I})$ the generalized madogram introduced in Fonseca et al. (2012), which in turns is the the $\lambda$-madogram (Naveu et al. (2009)) when $d=2=p$ and $1-\lambda_2=\lambda_1\in (0,1)$.\\ 
\noindent The max-min coefficient and the generalized madograms for pairs of sets $I_i$ and $I_j$ can be related throughout  
\begin{eqnarray*}\label{10}
R({\bf X},{\bf \lambda},\mathcal{I})\geq \displaystyle\bigvee _{\{i,j\}\subset \{1,...,p\}}R(({\bf X}_{I_i},{\bf X}_{I_j}),(\lambda_i,\lambda_j),\{I_i,I_j\}).
\end{eqnarray*}

\noindent We first present a key result that relates the expectation of $\displaystyle\bigvee _{j=1}^dF^{\lambda_j}(X_j)$ with the tail dependence function of $G$, which enables the derivation of the main properties of $R({\bf X},{\bf \lambda},\mathcal{I})$. The result also points out that in this work we can assume that the MEV distributions have unit Fr\'{e}chet margins without loss of generality.

\begin{pro}\label{p.1} Let ${\bf X}$ be a random vector with MEV distribution $G$, unit Fr\'{e}chet $F$ margins and tail dependence function $l$. If ${\bf Y}$ has MEV distribution with marginal distributions $F_j$, $j=1,...,d$, and the same copula as ${\bf X}$, then  for each ${\bf \lambda} \in (0,\infty)^d$, it holds that
\begin{eqnarray}\label{1}
E\left(\displaystyle\bigvee _{j=1}^dF_j^{\lambda_j}(Y_j)\right)=E\left(\displaystyle\bigvee _{j=1}^dF^{\lambda_j}(X_j)\right)=\frac{l(\lambda_1^{-1},...,\lambda_d^{-1})}{1+l(\lambda_1^{-1},...,\lambda_d^{-1})}.
\end{eqnarray}
\end{pro}

\dem We first deduce the distribution of $\displaystyle\bigvee _{j=1}^dF_j^{\lambda_j}(Y_j)$. Denoting the copula of ${\bf X}$ by $C_{{\bf X}}$,
we have, for each $u\in [0,1]$,
$$P\left(\displaystyle\bigvee _{j=1}^dF_j^{\lambda_j}(Y_j)\leq u\right)=C_{{\bf X}}\left(u^{\lambda_1^{-1}},...,u^{\lambda_d^{-1}}\right)=G\left(-\frac{1}{\lambda_1^{-1}\log u},...,\frac{1}{\lambda_d^{-1}\log u}\right)=$$

$$\exp \left(-l\left((-\log u)\lambda_1^{-1},...,(-\log u)\lambda_d^{-1}\right)\right)=u^{l\left(\lambda_1^{-1},...,\lambda_d^{-1}\right)}.$$

\noindent  Then 
$$E\left(\displaystyle\bigvee _{j=1}^dF_j^{\lambda_j}(Y_j)\right)=\int_0^1u^{l\left(\lambda_1^{-1},...,\lambda_d^{-1}\right)}l\left(\lambda_1^{-1},...,\lambda_d^{-1}\right)du=\frac{l(\lambda_1^{-1},...,\lambda_d^{-1})}{1+l(\lambda_1^{-1},...,\lambda_d^{-1})}.$$
\hspace{15cm}\fdem\\

\noindent The next result shows that the max-min coefficient takes into account the tail dependence function of all subvectors ${\bf X}_{\displaystyle\cup_{j\in T}I_j}$ of ${\bf X}$ with indices in $\displaystyle\cup_{j\in T}I_j$, $\emptyset \neq T\subseteq\{1,...,p\}$.\\

%$R({\bf X},{\bf \lambda})$, which corresponds to the particular case of $I_j=\{j\},\,j=1,...,p=d$.

\noindent From the Proposition  \ref{p.1} it holds that, for each $\emptyset \neq T\subseteq \{1,...,p\}$,\\
\begin{eqnarray}\label{8c}
e\left(I_j, j\in T\right)\equiv E\left(\displaystyle\bigvee _{j\in T}\displaystyle\bigvee _{i\in I_j}F^{\lambda_j}(X_i)\right)=
\frac{l\left(\displaystyle\sum_{j\in T}\lambda_j^{-1}\delta_1(I_j),...,\displaystyle\sum_{j\in T}\lambda_j^{-1}\delta_d(I_j)\right)}
{1+l\left(\displaystyle\sum_{j\in T}\lambda_j^{-1}\delta_1(I_j),...,\displaystyle\sum_{j\in T}\lambda_j^{-1}\delta_d(I_j)\right)},
\end{eqnarray}
\noindent leading to the following relations of the max-min coefficients with the tail dependence and the extremal coefficients.

\begin{pro}\label{p.4} If ${\bf X}$ has MEV distribution then, for each partition $\mathcal{I}=\{I_1,...,I_p\}$  of $D$ and ${\bf \lambda} \in (0,\infty)^p$, it holds that
\begin{eqnarray}\label{9}
R({\bf X},{\bf \lambda},\mathcal{I})=e\left(I_j, j\in \{1,...,p\}\right)-
\sum_{\emptyset \neq T\subseteq\{1,...,p\}}(-1)^{|T|+1}e\left(I_j, j\in T\right)
\end{eqnarray}

\noindent and

\begin{eqnarray}\label{9a}
 R({\bf X},{\bf 1},\mathcal{I})=\frac{\epsilon_{{\bf X}}}{1+\epsilon_{{\bf X}}}-
\sum_{\emptyset \neq T\subseteq\{1,...,p\}}(-1)^{|T|+1}\displaystyle\frac{\epsilon_{{\bf X}_{\displaystyle\cup_{j\in T}I_j}}}{1+\epsilon_{{\bf X}_{\displaystyle\cup_{j\in T}I_j}}}.
\end{eqnarray}
\end{pro}
\vspace{0.5cm}

\dem To obtain the first equality we first apply in (\ref{8b}) the relation 
\begin{eqnarray*}\label{3a}
\displaystyle\bigwedge _{j=1}^p\displaystyle\bigvee _{i\in I_j}F^{\lambda_j}(X_i)=\sum_{\emptyset \neq T\subseteq \{1,...,p\}}(-1)^{|T|+1}\displaystyle\bigvee _{j\in T}\displaystyle\bigvee _{i\in I_j}F^{\lambda_j}(X_i)
\end{eqnarray*}
\noindent  and then the Proposition 2.1 with (\ref{8c}). The  statement in (\ref{9a}) is a consequence of (\ref{9}) and (\ref{2a}).
\fdem\\

\noindent For the particular case of $I_j=\{j\},\,j=1,...,p=d$, we denote $R({\bf X},{\bf \lambda},\mathcal{I})$ simply by $R({\bf X},{\bf \lambda})$ and we have, as a consequence of the above result, that

\begin{eqnarray*}\label{2}
R({\bf X},{\bf \lambda})=\frac{l(\lambda_1^{-1},...,\lambda_d^{-1})}{1+l(\lambda_1^{-1},...,\lambda_d^{-1})}-
\sum_{\emptyset \neq S\subseteq\{1,...,d\}}(-1)^{|S|+1}\frac{l(\lambda_1^{-1}\delta_1(S),...,\lambda_d^{-1}\delta_d(S))}{1+l(\lambda_1^{-1}\delta_1(S),...,\lambda_d^{-1}\delta_d(S))}
\end{eqnarray*}
\noindent and
\begin{eqnarray*}\label{3}
 R({\bf X},{\bf 1})=\frac{\epsilon_{{\bf X}}}{1+\epsilon_{{\bf X}}}-
\sum_{\emptyset \neq S\subseteq\{1,...,d\}}(-1)^{|S|+1}\displaystyle\frac{\epsilon_{{\bf X}_S}}{1+\epsilon_{{\bf X}_S}}.
\end{eqnarray*}\\

\noindent These two relations extend the result of the Proposition 1 in Naveau et al. (2009) and equation (14) in Cooley et al. (2006), where, for $d=2$, $1-\lambda_2=\lambda_1=\lambda \in (0,1)$, we have 

\begin{eqnarray*}\label{5}
\frac{1}{2}R({\bf X},{\bf \lambda})=\frac{l(\frac{1}{\lambda},\frac{1}{1-\lambda})}{1+l(\frac{1}{\lambda},\frac{1}{1-\lambda})}-\frac{3}{2}\frac{1}{(1+\lambda)(2-\lambda)}
\end{eqnarray*}\\

\noindent and $R({{\bf X}},{\bf 1})=\frac{\epsilon_{{\bf X}}-1}{\epsilon_{{\bf X}}+1}$.\\

\noindent Our next step is to compare the value of the max-min coefficient $R({\bf X},{\bf \lambda},\mathcal{I})$ with the co\-rres\-ponding coefficient in the two boundary cases of  independent or totally dependent ${\bf X}_{I_j}$, $j=1,...,p$.\\

%$$$$$$$$$$
\begin{pro}\label{p.5} Let  ${\bf X}$, ${\hat{{\bf X}}}=(\hat{X}_1,...,\hat{X}_d)$ and ${\bar{{\bf X}}}=(\bar{X}_1,...,\bar{X}_d)$ be vectors of unit Fr\'{e}chet random variables with  MEV distributions such that $\hat{{\bf X}}_{I_j}$, $j=1,...,p$, are independent, $\bar{{\bf X}}_{I_j}$, $j=1,...,p$, are totally dependent and, for each $j=1,...,p$, $\hat{{\bf X}}_{I_j}$, $\bar{{\bf X}}_{I_j}$ and  ${\bf X}_{I_j}$ are identically distributed. Then, for each ${\bf \lambda} \in (0,\infty)^d$ and partition $\mathcal{I}$ of $D$, it holds that\\

\noindent (a) $R({\bar{{\bf X}}},{\bf \lambda},\mathcal{I})\leq R({\bf X},{\bf \lambda},\mathcal{I})\leq R({\hat{{\bf X}}},{\bf \lambda},\mathcal{I})$,\\

\noindent (b) $R({\bar{{\bf X}}},{\bf 1},\mathcal{I})=0,$\\

\noindent (c) $R({\hat{{\bf X}}},{\bf 1},\mathcal{I})=\frac{\displaystyle\sum_{j=1}^p\epsilon_{{\bf X}_{I_j}}}{1+\displaystyle\sum_{j=1}^p\epsilon_{{\bf X}_{I_j}}}-\displaystyle\sum_{\emptyset \neq T\subseteq\{1,...,p\}}(-1)^{|T|+1}\frac{\displaystyle\sum_{j\in T}\epsilon_{{\bf X}_{I_j}}}{1+\displaystyle\sum_{j\in T}\epsilon_{{\bf X}_{I_j}}}.$
\end{pro}

\dem If ${\bf X}$  has MEV distribution then it is associated (Marshall and Olkin (1983)) and then the variables $M(I_j)$, $j=1,...,p$, are also associated (Esary et al. (1967)). For these associated variables it holds that

$$\displaystyle\bigwedge_{j=1}^pP\left(M(I_j)>x_j\right)\geq 
P\left(\displaystyle\bigcap_{j=1}^p\{M(I_j)>x_j\}\right)\geq \displaystyle\prod_{j=1}^pP\left(M(I_j)>x_j\right)$$

\noindent and

$$\displaystyle\bigwedge_{j=1}^pP\left(M(I_j)\leq x_j\right)\geq P\left(\displaystyle\bigcap_{j=1}^p\{M(I_j)\leq x_j\}\right)\geq \displaystyle\prod_{j=1}^pP\left(M(I_j)\leq x_j\right).$$
\vspace{0.3cm}

\noindent By taking ${\hat{M}}(I_j)=\displaystyle\bigvee_{i\in I_j}{\hat{X}}_i$ and ${\bar{M}}(I_j)=\displaystyle\bigvee_{i\in I_j}{\bar{X}}_i$, $j=1,...,p$, we can rewrite the above inequalities as 

$$P\left(\displaystyle\bigcap_{j=1}^p\{{\bar{M}}(I_j)>x_j\}\right) \geq P\left(\displaystyle\bigcap_{j=1}^p\{M(I_j)>x_j\}\right) \geq P\left(\displaystyle\bigcap_{j=1}^p\{{\hat{M}}(I_j)>x_j\}\right)$$

\noindent and

$$P\left(\displaystyle\bigcap_{j=1}^p\{{\bar{M}}(I_j)\leq x_j\}\right) \geq 
P\left(\displaystyle\bigcap_{j=1}^p\{M(I_j)\leq x_j\}\right) \geq P\left(\displaystyle\bigcap_{j=1}^p\{{\hat{M}}(I_j)\leq x_j\}\right)$$

\noindent This concordance order implies (Shaked and Shanthikumar (2007)) that 
$$E\left(f\left(\displaystyle\bigvee_{j=1}^p{\bar{M}}(I_j)\right)\right)\leq E\left(f\left(\displaystyle\bigvee_{j=1}^pM(I_j)\right)\right)
\leq E\left(f\left(\displaystyle\bigvee_{j=1}^p{\hat{M}}(I_j)\right)\right),$$
\noindent and 
$$E\left(f\left(\displaystyle\bigwedge_{j=1}^p{\bar{M}}(I_j)\right)\right)\geq E\left(f\left(\displaystyle\bigwedge_{j=1}^pM(I_j)\right)\right)
\geq E\left(f\left(\displaystyle\bigwedge_{j=1}^p{\hat{M}}(I_j)\right)\right)$$

\vspace{0,3cm}

\noindent for all non-decreasing functions $f$. The result in (a) then follows by taking $f=F$ and replacing $X_i$ by  $\frac{X_i}{\lambda_j}$, $i\in I_j$, for each $j=1,...,p$. 

\noindent The equalities in (b) and (c) follows from (\ref{9a}) and, in particular for (b), we recall that if   $\bar{{\bf X}}_{I_j}$, $j=1,...,p$, are totally dependent vectors then the copula of ${\bar{{\bf X}}}$ is also the copula of the minimum (Nelsen, 2006). \fdem\\

\noindent For the particular case of $I_j=\{j\}$, $j=1,...,d$, the equality in (c) leads to 
$$R({\hat{{\bf X}}},{\bf 1})=\frac{d}{d+1}-\displaystyle\sum_{\emptyset \neq S\subseteq\{1,...,d\}}(-1)^{|S|+1}\frac{|S|}{1+|S|}
=\frac{d}{d+1}-\displaystyle\sum_{k=1}^d(-1)^{k+1}\left(_k^d\right)\frac{k}{k+1}=\frac{d-1}{d+1}.$$

%$$\frac{d}{d+1}+\frac{1}{d+1}\displaystyle\sum_{k=1}^d(-1)^{k+2}\left(_{k+1}^{d+1}\right)(k+1)-\frac{1}{d+1}\displaystyle\sum_{k=1}^d(-1)^{k+2}\left(_{k+1}^{d+1}\right)=$$

%$$\frac{d}{d+1}+\frac{1}{d+1}\left(0-\left(_{1}^{d+1}\right)\right)+\frac{1}{d+1}\left(0-\left(_{0}^{d+1}\right)+\left(_{1}^{d+1}\right)\right)=\frac{d}{d+1}+\frac{1}{d+1}\displaystyle\sum_{k=1}^d(-1)^{k+2}\left(_{k+1}^{d+1}\right)k=\frac{d-1}{d+1}.$$

\noindent which extends the result for the case of $d=2$, where $R({\hat{{\bf X}}},{\bf 1})=\frac{1}{3}$.\\

%%%%%%%%%%%%%%%%%%%%%%%%%%%%%%%%%%%%%%%%%%%%%%%%%%%%%%%%%%%%%%%%%%%%
\section{Examples and an application} 
 In order to illustrate the previous results, we consider two families of  Multivariate Extreme Value distributions and we compute the expressions for $R({\bf X},{\bf \lambda},\mathcal{I})$, which can be easily implemented.

\begin{ex}  If {\bf X} is the MEV marginal distribution of the Multivariate Maxima of Moving Maxima processes considered in Smith and Wissman (1996) then $l(x_1,...,x_d)=\displaystyle\sum_{l=1}^\infty \sum_{k=-\infty}^\infty\bigvee_{j=1}^d x_j\,\alpha_{l,k,j}$. For this tail dependence function we obtain\\

\noindent $R({\bf X},{\bf \lambda},\mathcal{I})=\displaystyle\sum_{\emptyset \neq T\subsetneq\{1,...,p\}}\frac{(-1)^{|T|+1}}{1+
\displaystyle\sum_{l=1}^\infty \displaystyle\sum_{k=-\infty}^\infty\bigvee_{t\in T} \bigvee_{j\in I_t}\lambda_t^{-1}\,\alpha_{l,k,j}}
-\frac{1+(-1)^p}{1+
\displaystyle\sum_{l=1}^\infty \displaystyle\sum_{k=-\infty}^\infty\bigvee_{t=1}^p \bigvee_{j\in I_t}\lambda_t^{-1}\,\alpha_{l,k,j}}$\\

\noindent and, in particular,\\

\noindent $R({\bf X},{\bf \lambda})=\displaystyle\sum_{\emptyset \neq S\subsetneq\{1,...,d\}}\frac{(-1)^{|S|+1}}{1+
\displaystyle\sum_{l=1}^\infty \displaystyle\sum_{k=-\infty}^\infty\bigvee_{j\in S} \lambda_j^{-1}\,\alpha_{l,k,j}}
-\frac{1+(-1)^d}{1+
\displaystyle\sum_{l=1}^\infty \displaystyle\sum_{k=-\infty}^\infty\bigvee_{j=1}^d \lambda_j^{-1}\,\alpha_{l,k,j}}.$ \fdem

\end{ex}

\begin{ex} For the Symmetric Logistic model $l(x_1,...,x_d)=\left(\displaystyle\sum_{j=1}^dx_j^{1/\theta}\right)^\theta$ we have\\

\noindent $R({\bf X},{\bf \lambda},\mathcal{I})=\displaystyle\sum_{\emptyset \neq T\subsetneq\{1,...,p\}}\frac{(-1)^{|T|+1}}{1+
\left(\displaystyle\sum_{t\in T}\lambda_t^{-1/\theta}|I_t|\right)^\theta}
-\frac{1+(-1)^p}{1+
\left(\displaystyle\sum_{t=1}^p\lambda_j^{-1/\theta}|I_t|\right)^\theta},$\\
\vspace{0.7cm}

\noindent $R({\bf X},{\bf \lambda})=\displaystyle\sum_{\emptyset \neq S\subsetneq\{1,...,d\}}\frac{(-1)^{|S|+1}}{1+
\left(\displaystyle\sum_{j\in S}\lambda_j^{-1/\theta}\right)^\theta}
-\frac{1+(-1)^d}{1+
\left(\displaystyle\sum_{j=1}^d\lambda_j^{-1/\theta}\right)^\theta}$\\

\noindent and\\

\noindent $R({\bf X},{\bf 1})=\displaystyle\sum_{k=1}^{d-1}(-1)^{k+1}\left(^d_k\right)\frac{1}{1+k^\theta}
-(1+(-1)^d)\frac{1}{1+
d^\theta}$. \fdem

\end{ex} 
\vspace{0.5cm}

%__________________

\noindent Several parametric and non-parametric estimators for the tail dependence function are available in the literature (Beirlant et al. (2004), Schmidt and Stadtmuller (2006), Krajina (2010)) which can be applied to the terms in (\ref{9}). The comparison of estimation procedures is out of our purposes in this paper and we simply remark that the definition of the max-min dependence coefficient suggests a non-parametric estimator based on sample means.\\

\noindent Let ${\bf X}^{(k)}=(X_1^{(k)},...,X_d^{(k)})$, $k=1,...,n$, be a sequence of independent copies of ${\bf X}$ and ${\hat F}_j$  the empirical distribution provided by $X_j^{(k)}$, $k=1,...,n$, $j=1,...,d$.\\
A natural estimator for $R({\bf X},{\bf \lambda},\mathcal{I})$ is 
$${\hat R}({\bf X},{\bf \lambda},\mathcal{I})=\frac{1}{n}\displaystyle\sum_{k=1}^n\left(\displaystyle\bigvee _{j=1}^p\displaystyle\bigvee _{i\in I_j}{\hat F}_i^{\lambda_j}(X_i^{(k)})-\displaystyle\bigwedge _{j=1}^p\displaystyle\bigvee _{i\in I_j}{\hat F}_i^{\lambda_j}(X_i^{(k)})\right),$$

\noindent and, in particular for $p=d$, 

$${\hat R}({\bf X},{\bf \lambda})=\frac{1}{n}\displaystyle\sum_{k=1}^n\left(\displaystyle\bigvee _{j=1}^d{\hat F}_j^{\lambda_j}(X_j^{(k)})-\displaystyle\bigwedge _{j=1}^d{\hat F}_j^{\lambda_j}(X_j^{(k)})\right).$$

\noindent If we denote ${\hat M}_k(I_j)^{\lambda_j}=\displaystyle\bigvee _{i\in I_j}{\hat F}_i^{\lambda_j}(X_i^{(k)})$ 
 then we can write 
 \begin{eqnarray*}\label{11}
{\hat R}({\bf X},{\bf \lambda},\mathcal{I})=
\frac{1}{n}\displaystyle\sum_{k=1}^n\displaystyle\bigvee _{j=1}^p{\hat M}_k(I_j)^{\lambda_j}
-\sum_{\emptyset \neq T\subseteq\{1,...,p\}}(-1)^{|T|+1}\displaystyle\frac{1}{n}\sum_{k=1}^n
\displaystyle\bigvee _{j\in  T}{\hat M}_k(I_j)^{\lambda_j}.
\end{eqnarray*}

\noindent The strong consistency of the terms of this sum is stated in the proof of the Proposition 3.8  of Ferreira and Ferreira (2012b) and the asymptotic normality can be deduced from the Theorem 6 in Fermanian et al. (2004).\\

\noindent As an application of this estimation procedure we consider for  ${\bf X}$ some financial stock markets grouped in the tree big world markets $I_1=$Europe, $I_2=$USA and $I_3=$Far East, as considered in Ferreira and Ferreira (2012b). The data are monthly maximums of the negative log-returns of the closing values of the stock market indexes CAC 40 (France), FTSE100 (UK), SMI (Swiss), XDAX (German), Dow Jones (USA), Nasdaq (USA), SP500 (USA), HSI (China) and Nikkei (Japan), from January 1993 to March 2004.\\

\noindent  In the table 1 we present the estimates corresponding to  ${\bar M}(A)=\frac{1}{n}\displaystyle\sum_{k=1}^n\displaystyle\bigvee _{i \in A}{\hat F}_i(X_j^{(k)})$  which we need to compute ${\hat R}({\bf X},{\bf 1},\mathcal{I})$ for $\mathcal{I}=\{I_1,I_2,I_3\}$.\\
 
\noindent  We obtain ${\hat R}({\bf X},{\bf 1},\mathcal{I})=0,321$, ${\hat R}(({\bf X}_{I_1},{\bf X}_{I_2}),{\bf 1},\{I_1,I_2\})=0,172$, ${\hat R}(({\bf X}_{I_1},{\bf X}_{I_3}),{\bf 1},\{I_1,I_3\})=0,222$ and ${\hat R}(({\bf X}_{I_2},{\bf X}_{I_3}),{\bf 1},\{I_2,I_3\})=0,247$.\\

\begin{table}[h]
\centering
\footnotesize
\begin{tabular}{| c | c | c | c | c | c | c | c  |}
\hline
$A$ & $I_1$ & $I_2$ & $I_3$ & $I_1 \cup I_2$ & $I_1 \cup I_3$ & $I_2 \cup I_3$ & $I_1 \cup I_2  \cup I_3$ \\  \hline
${\bar M}(A)$ & 0,691736695    &   0,614005602     &  0,625910364 & 0,738655462 & 0,770028011 & 0,743557423 &  0,80070028\\ \hline
\end{tabular}
\caption{$I_1=$Europe, $I_2=$USA and $I_3=$Far East}
\end{table}

%%%%%%%%%%%%%%%%%%%%%%%%%%%%%%%%%%%%%%%%%%%%%%%%%%%%%%%%%%%%%%%%
%%%%%%%%%%%%%%%%%%%%%%%%%%%%%%%%%%%%%%%%%%%%%%%%%%%%%%%%%%%%%%%

\end{document}